\documentclass{article}
\usepackage{amssymb}

\newcommand {\Q}{{\mathbb{Q}}}

\newcommand {\R}{{\mathbb{R}}}
\newcommand {\Z}{{\mathbb{Z}}}
\newcommand {\N}{{\mathbb{N}}}

\newtheorem{theorem}{Theorem}
\newtheorem{definition}{Definition}

\newtheorem{lemma}{Lemma}
\newcommand{\eindebewijs}{\hfill$\Box$\par\medskip}
\begin{document}
\title{Solving conics over $\Q(t_1,\ldots,t_k)$.}
\author{Mark van Hoeij \\
Department of Mathematics \\ 
Florida State University \\
Tallahassee, FL 32306-3027, USA \\
hoeij@math.fsu.edu}
\maketitle
\begin{abstract}
Let $K = \Q(t_1,\ldots,t_k)$ and $a,b,c \in K$. We give a simple
algorithm to find, if it exists, $X,Y,Z$ in $K$, not all zero,
for which $aX^2 + bY^2 + cZ^2 = 0$.
\end{abstract}

\section{Introduction}
\label{section1}
Let $K = \Q(t_1,\ldots,t_k)$ and $a,b,c \in K$. Our goal is to
find, if it exists, a solution $X,Y,Z \in K$, not all zero, for which
\begin{equation}
	\label{conic}
	aX^2 + bY^2 + cZ^2 = 0.
\end{equation}
We will use projective coordinates $(X:Y:Z) \in P^2(K)$ to indicate that
at least one of $X,Y,Z \in K$ must be non-zero,
and that $(sX: sY: sZ)$ will be identified
with $(X:Y:Z)$ for any nonzero $s \in K$.

The case $k=0$ ($K=\Q$) was solved by Legendre, who gave a criterion
(local-global principle for conics) to decide if a solution
$(X,Y,Z) \in P^2(\Q)$ exists. He also gave an algorithm to
find a solution. 

Legendre showed that a solution $(X:Y:Z) \in P^2(\Q)$ exists if
and only if a so-called {\em solubility certificate} exists.
The non-trivial part is to show that existence of this
certificate implies the existence of a solution $(X:Y:Z)$.
We sketch a textbook proof, with an improvement given in \cite{cm},
because our algorithm will follow the same approach.
Use the solubility certificate to construct a certain lattice
$L$ which has the property that $aX^2+bY^2+cZ^2 \equiv 0$ mod $abc$ for
every $(X,Y,Z) \in L$, and use Minkowski's theorem to show that
$L$ contains a non-zero element $(X,Y,Z)$
for which $|a|X^2 + |b|Y^2 + |c|Z^2 < 2abc$. If $aX^2 +bY^2 + cZ^2$ is
not already 0, it must be $\pm abc$, and a transformation yields new $X,Y,Z$
for which $aX^2+bY^2+cZ^2=0$. Alternatively, in \cite{cm} it is
shown that one can construct a lattice $L' \subset L$ such
that $aX^2+bY^2+cZ^2 \equiv 0$ mod $2abc$ for every $(X,Y,Z) \in L'$
and moreover that this $L'$ still contains a non-zero element $(X,Y,Z)$
for which $|a|X^2 + |b|Y^2 + |c|Z^2 < 2abc$ and hence $aX^2+bY^2+cZ^2=0$.
Given a solubility certificate, one can quickly
find such $(X,Y,Z)$ with the \cite{LLL} lattice reduction algorithm
and Lemma 2.7 in \cite{cremona} (where an even faster algorithm
was given as well).

In the proof of Theorem~\ref{thm1} we will follow the same ideas for
$K = \Q(t_1,\ldots,t_k)$. Although this involves additional notation,
some details will be easier since the lattice arguments for $k=0$ simplify
to linear algebra when $k>0$.

\section{The solubility certificate}
\label{cert}
If $abc=0$ then the problem is trivial (if for example $a=0$ then
$(1:0:0)$ is a solution).
If $k=0$ and if $a,b,c$ have the same sign then there is no
solution, not even over $\R$. From now on we assume
\begin{equation}
	\label{nonz}
	abc \neq 0, \
  \mbox{\rm and if}\ k=0 \ \mbox{\rm then}\ a,b,c \ \mbox{\rm do not all have the same sign}.
\end{equation}
After multiplying by denominators we may assume that
\begin{equation}
	\label{Rk}
	a,b,c \in R_k := \left\{ \begin{array}{ll}
\Z & {\rm if \ } k=0 \\
\Q[t_1,\ldots,t_k] & {\rm if \ } k>0
 \end{array}  \right.
\end{equation}
Let $g = {\rm gcd}(a,b)$, where gcd refers to the greatest common
divisor in $R_k$.
If $g \neq 1$ then $a,b,c$ will be replaced by $a/g, b/g, cg$.
The new conic is equivalent to the old conic under
$(X:Y:Z) \mapsto (X:Y:Z/g)$.
The same is done for ${\rm gcd}(b,c)$ and ${\rm gcd}(c,a)$,
after which we get
\begin{equation}
	\label{gcd}
	{\rm gcd}(a,b) = {\rm gcd}(b,c) = {\rm gcd}(c,a) = 1
\end{equation}

A non-trivial square is a square of a non-zero $d \in R_k$ that is not a unit
(so $d^2 \neq 1$ if $k=0$ and $d \not\in \Q$ if $k>0$).
If $a$, $b$ or $c$ is divisible by a non-trivial square
then we can divide and obtain an equivalent conic.
Thus, we may assume that $a$, $b$, $c$ are not divisible by
non-trivial squares, which because of (\ref{gcd}) is
equivalent to
\begin{equation}
	\label{square}
	abc \, \ \mbox{{\rm is not divisible by a non-trivial square.}}
\end{equation}
Given equation~(\ref{conic}), the above describes Legendre's procedure
to compute a {\em reduced form} of this equation, which means:
\begin{definition}
\label{reduced}
An equation $aX^2 + bY^2 + cZ^2 = 0$ is in {\em reduced form} if $a,b,c$
satisfy assumptions (\ref{nonz}), (\ref{Rk}), (\ref{gcd}) and (\ref{square}).
\end{definition}

In this paper, if $k>0$, then the {\em leading coefficient} resp. {\em degree}
of a polynomial refers to the leading coefficient resp. degree
with respect to $t_k$.
If $k>0$ then denote $d_a,d_b,d_c \in \N$ as the degrees of $a,b,c$.
Denote $l_a, l_b, l_c \in \Q[t_1,\ldots,t_{k-1}]$ as the leading coefficients of $a,b,c$.
Set
\[
	{\rm case} :=
	\left\{ \begin{array}{ll}
  0 & \mbox{\rm if}\ k>0 \ \mbox{\rm and}\
  d_a \equiv d_b \equiv d_c {\rm \ mod \ }2 \\
  1 & {\rm otherwise} \end{array} \right.
\]
If a solution $(X:Y:Z) \in P^2(K)$ of equation~(\ref{conic}) exists, then
\begin{equation}
	\label{lc}
	\left\{ \begin{array}{l}
	{\rm case} = 1, \ {\bf \rm or} \\
	l_a x^2 + l_b y^2 + l_c z^2 = 0 \ \mbox{\rm has a solution in}\ P^2( \Q(t_1,\ldots,t_{k-1}) ).
	\end{array} \right.
\end{equation}
{\bf Proof:} Assume $k>0$. After scaling we may assume that $X,Y,Z \in R_k$.
Let $d$ be the maximum of the degrees (always with respect to $t_k$)
of $aX^2$, $bY^2$, and $cZ^2$. If the degree of $aX^2$ is $d$ then let $x$ be
the leading coefficient of $X$, otherwise set $x=0$. By doing the same for $y$ and $z$
one finds a solution $(x:y:z)$. \eindebewijs
\vspace{10pt}
\noindent Denote
\[
	\tilde{R}_k := \left\{ \begin{array}{ll}
\Z & {\rm if \ } k=0 \\
\Q(t_1,\ldots,t_{k-1})[t_k] & {\rm if \ } k>0
 \end{array}  \right.
\]
Let ${\rm Factors}(a)$ denote the set of all irreducible $p \in \tilde{R}_k$
that divide $a$, with the additional condition that $p$ must be positive
(and hence a prime number) when $k=0$,
and $p$ must be monic (leading coefficient 1) when $k>0$. Denote
\[ {\bf F}_p := \tilde{R}_k/(p). \]
Introduce a new variable $T$ and let $f_a,f_b,f_c \in R_k[T]$ be
the following polynomials
\[
        f_a := bT^2 + c, \ \ \ \ f_b := cT^2 + a,  \ \ \ \ f_c := aT^2 + b.
\]
For $f \in \tilde{R}_k[T]$, the image in ${\bf F}_p[T]$ is denoted as $f$ mod $p$.

\begin{definition}
\label{defcert}
Let $aX^2 + bY^2 + cZ^2 = 0$ be in reduced form (see Definition~\ref{reduced}).
Then a {\em solubility certificate} is
a list containing the following: \\[3pt]
\mbox{} \ \  For every $p \in {\rm Factors}( \hspace{-0.5pt} a \hspace{-0.5pt} )$,
	a root of $f_a$\hspace{-0.8pt} mod $p$ in ${\bf F}_p$ \\
\mbox{} \ \  For every $p \in {\rm Factors}(b)$, a root of $f_b$ mod $p$ in ${\bf F}_p$ \\
\mbox{} \ \  For every $p \in {\rm Factors}(c)$, a root of $f_c$ mod $p$ in ${\bf F}_p$ \\
\mbox{} \ \  If\, ${\rm case}=0$, a reduced form plus solubility certificate for equation~(\ref{lc}).
\end{definition}
Note that for $k=0$, the standard
definition is formulated differently, namely the roots
are combined by Chinese remaindering. So the standard definition of a solubility certificate
is a list of three integers $k_a, k_b, k_c$ that satisfy $f_a$ mod $a$, $f_b$ mod $b$,
and $f_c$ mod $c$ respectively.

\begin{lemma}
Assume equation~(\ref{conic}) is in reduced form and has a
solution in $P^2(K)$. Then a solubility certificate exists.
\end{lemma}
{\bf Proof:}
Let $(X:Y:Z)$ be a solution. After scaling
we may assume that $X,Y,Z \in R_k$ with ${\rm gcd}(X,Y,Z)=1$.
For ${\rm case}=0$ we have already shown that equation~(\ref{lc}) has a solution by
considering the leading coefficients of $X,Y,Z$.
So equation~(\ref{lc}) has a solubility certificate by induction.
The remainder of the proof is the same as for the $k=0$ case,
we include it for completeness since our notations are slightly different.

Let $p \in {\rm Factors}(a)$.
Now $0 = aX^2 + bY^2 + cZ^2$ reduces to $bY^2 + cZ^2$ mod $p$.
Since $b,c$ do not vanish mod $p$, it follows that
$\overline{Y}, \overline{Z}$
are either both zero or both not zero, where $\overline{Y}, \overline{Z}$
denote the images of $Y,Z$ in ${\bf F}_p$.
The first case is ruled out since it would imply that $bY^2 + cZ^2$ and hence $aX^2$
would be divisible by $p^2$, which implies that $X$ is divisible by $p$,
contradicting ${\rm gcd}(X,Y,Z)=1$. Thus, $\overline{Y}, \overline{Z}$ are not zero.
Then $f_a$ mod $p$ has $\pm \overline{Y}/\overline{Z} \in {\bf F}_p$ as roots.
Repeating this for $f_b$ resp. $f_c$ for the factors of $b$ resp. $c$
it follows that a solubility certificate exists.
\eindebewijs

\subsection{Algorithmic remarks}
\label{sect21}
To compute (if it exists) a solubility certificate we use:
\begin{enumerate}
\vspace{-2pt}
\item A factorization algorithm for $\tilde{R}_k$.
\vspace{-2pt}
\item A factorization algorithm for ${\bf F}_p[T]$.
\vspace{-2pt}
\item If ${\rm case}=0$ we need items 1 and 2 for $k-1$ as well.
\vspace{-2pt}
\end{enumerate}
Such algorithms are available in several computer algebra systems.
If $k=0$ then item 1 can take much time,
while the time spent on item 2 is negligible in comparison.
If $k>0$, it is more or less the other way around. When $k>0$ the time
spent on items 1 and 2 is polynomially bounded in terms of the
degrees and logarithmic heights, provided that $k$ is fixed.
So for fixed $k$, one can compute a solubility certificate
in polynomial time {\em except}\,  if ${\rm case}=0$ at each of the stages
$k,k-1,\ldots,2,1$ because then we need to factor in $\tilde{R}_0 = \Z$
for which no polynomial time algorithm is known.
Having to factor integers
can sometimes, but not always, be avoided by interchanging variables. \\

\noindent{\bf Remark. A solubility certificate corresponds to a subset of the conic.} \\
If equation~(\ref{conic}) has a solution in $P^2(K)$ then the number
of solubility certificates equals $2^n$ for some integer $n$,
because for each $p$, one has to choose one of the two roots (one root if $p=2$).
Every certificate corresponds
to a non-empty (see Theorem~\ref{thm1}) subset of the conic, where the word
conic refers to the set of solutions of equation~(\ref{conic}) in $P^2(K)$.
Grouping together those subsets/certificates that correspond under $(X:Y:Z)
\mapsto (\pm X : \pm Y : \pm Z)$, one writes the conic as a disjoint
union of non-empty subsets $S_1, S_2,\ldots,S_{2^m}$ for some integer $m \geq 0$.
Given one point on the conic, one can write down a {\em parametrization}, a
birational map $\Psi$ from $P^1(K)$ to the conic. Then every point on the conic
can be written as $\Psi(u:v)$ for some $(u:v) \in P^1(K)$. Nevertheless, one can
construct examples with $k=0$, $m>0$, with a point in one $S_i$ given, where it would
require a breakthrough in integer factorization to find a point in any $S_j$ with $j \neq i$.
This is because going from one $S_i$ to another is equivalent to finding
a non-trivial factor of $a$, $b$ or $c$. \\


Although the above remark is not needed for our algorithm, it is useful to know
that a solubility certificate corresponds to a subset of the conic, and that the
algorithm in the next section finds an element of that subset. So
any certificate will lead to a point on the conic, however,
a much nicer point (smaller degrees
or heights) might be found by trying $2^m$ certificates.

\section{Finding a point on the conic}
Let $K = \Q(t_1,\ldots,t_k)$ where $t_1,\ldots,t_k$ are variables and $k \in \N$. \\[10pt]
{\nopagebreak
{\bf Algorithm FindPoint} \\[2pt]
{\bf Input:} $k$, \ $t_1,\ldots,t_k$, \ $a,b,c \in K$ satisfying assumptions
(\ref{nonz},\ref{Rk},\ref{gcd},\ref{square}), \ and a \\ 
\mbox{} \hspace{29.5pt} solubility certificate. \\[2pt]
{\bf Output:} A solution $(X:Y:Z) \in P^2(K)$ of $aX^2 + bY^2 + cZ^2 = 0$. }
\begin{enumerate}
\item If $k=0$ then use Algorithm II or III in \cite{cremona} and stop.
\item Let $d_a$, $d_b$, $d_c$ be the degrees (with respect to $t_k$) of $a$, $b$, $c$.
\item If $d_a \equiv d_b \equiv d_c$ mod 2 then ${\rm case} := 0$, otherwise ${\rm case} := 1$ and $E := \emptyset$.
\item Let $A$, $B$, $C$ be:
\[ \lceil \frac{d_b + d_c}2 \rceil - {\rm case}, \ \ \ \ 
   \lceil \frac{d_c + d_a}2 \rceil - {\rm case}, \ \ \ \
   \lceil \frac{d_a + d_b}2 \rceil - {\rm case}
\]
and set
\[  X := \sum_{i=0}^A X_i t_k^i, \ \ \ \ \ \ \
    Y := \sum_{i=0}^B Y_i t_k^i, \ \ \ \ \ \ \ \ \,
    Z := \sum_{i=0}^C Z_i t_k^i \,\ \ \ \mbox{} \]
where the $X_i,Y_i,Z_i$ are new variables.
\item \label{step5} If ${\rm case}=0$, let $l_a, l_b, l_c \in \Q[t_1,\ldots,t_{k-1}]$
be the leading coefficients of $a,b,c$. Let $(x:y:z)$ be a solution of
equation~(\ref{lc}), see the Remark below.
Introduce a new variable $w$ and set $E := \{X_A - xw, Y_B -  yw, Z_C - zw\}$.
\item \label{step6}
For every $p \in {\rm Factors}(a)$, let $\alpha \in {\bf F}_p$ be a root of $f_a$ mod
$p$. Lift $\alpha$ to $\tilde{R}_k$, divide $Y-\alpha Z$ by $p$ and let $r$
be the remainder. Write $r = \sum_{i<d} r_i t_k^i$ where $d$ is the degree of $p$.
Then $r_i$ is a $\Q(t_1,\ldots,t_{k-1})$-linear combination of $Y_j$'s and $Z_j$'s. Now
set $E := E \bigcup \{ r_0, \ldots, r_{d-1} \}$. \\
Repeat this for $f_b$ resp. $f_c$ for all $p$ in Factors$(b)$ resp. Factors$(c)$.
\item \label{step7} Equate all elements of $E$ to 0, resulting in a system of homogeneous
linear equations for the variables $X_i,Y_i,Z_i$ (and $w$ if\, ${\rm case} = 0$).
Solve this system and substitute the general solution into $X,Y,Z$. Multiply $(X,Y,Z)$ by the
least common multiple of the denominators of $X,Y,Z$.
Then divide by ${\rm gcd}(X,Y,Z)$. Return $(X:Y:Z)$ and stop.
\end{enumerate}

\noindent {\bf Remark on step~\ref{step5}:}
For the algorithm it is more convenient to replace
``a reduced form and a solubility certificate for equation~(\ref{lc})''
in  Definition~\ref{defcert}
by ``a solution $(x:y:z)$ of~(\ref{lc})''. Induction resp. recursion shows
that this is mathematically resp. algorithmically equivalent.

\begin{theorem}
\label{thm1}
Assume equation~(\ref{conic}) is in reduced form. Then it has a solution in $P^2(K)$
if and only if a solubility certificate exists.
\end{theorem}
{\bf Proof:} If $k=0$ then this is Legendre's theorem.
Assume we have a solubility certificate and that $k>0$. It suffices to prove that
Algorithm FindPoint works. For step~\ref{step5} see the Remark above.
The general solution of system $E$ in
step~\ref{step7} is not zero because the number of variables is the number of
equations plus 1. Therefore, the notation $(X:Y:Z)$ in step~\ref{step7} is
correct, that is, $X,Y,Z$ are not all zero.
Let $D$ be the degree of $abc$.
By construction, $aX^2 + bY^2 + cZ^2$ has degree $\leq {\rm max}\{d_a + 2A,
d_b + 2B, d_c + 2C\} = D-{\rm case}$.
Furthermore, it is divisible by every $p$ in ${\rm Factors}(a)$,
${\rm Factors}(b)$ and ${\rm Factors}(c)$ because equations $E$ are satisfied.
Then assumption~(\ref{square}) from section~\ref{cert} implies that
$aX^2 + bY^2 + cZ^2$ is divisible by $abc$ hence it can be written
as $s \cdot abc$ for some $s$. Comparing degrees, it follows that $s=0$
when ${\rm case}=1$. If ${\rm case}=0$, then the $t_k^D$ coefficient
in $aX^2 + bY^2 + cZ^2$ is (see step~\ref{step5}) $w^2$ times
the left-hand side of equation~(\ref{lc}), which equals 0.
It follows that $aX^2 + bY^2 + cZ^2 = 0$.

The free variable(s) in the general solution of $E$ will disappear
in step~\ref{step7} when we divide by
gcd$(X,Y,Z)$ because the conic does not contain any lines
(it is non-degenerate by assumption~(\ref{nonz})).
So the solution of equation~(\ref{conic}) returned by the algorithm
will be an element of $P^2(K)$.
\eindebewijs

\section{Algorithmic details}
Let $K = \Q(t_1,\ldots,t_k)$ where $t_1,\ldots,t_k$ are variables and $k \in \N$. \\[10pt]
{\nopagebreak
{\bf Algorithm ReduceConic} \\[2pt]
{\bf Input:} $k$, \ $t_1,\ldots,t_k$, \ \ $a,b,c \in K$ satisfying
assumption~(\ref{nonz}). \\[2pt]
{\bf Output:} $a_1,b_1,c_1$ for which $a_1X^2 + b_1Y^2 + c_1 Z^2=0$ is the reduced form of~(\ref{conic}).
\begin{enumerate}
\item Multiply $(a,b,c)$ by the least common multiple of the denominators.
\item $g := {\rm gcd}(a,b)$.  Set $a := a/g$, $b := b/g$, $c := cg$, and store $g$.
\item $g := {\rm gcd}(b,c)$.\hspace{1pt}  Set $b := b/g$, $c := c/g$, $a := ag$, and store $g$.
\item $g := {\rm gcd}(c,a)$.  Set $c := c/g$, $a := a/g$, $b := bg$, and store $g$.
\item If $k=0$ then let $A,B,C$ be the absolute values of $a,b,c$ and do
\begin{enumerate}
\item Compute all prime factors of $A$ that can be calculated reasonably quickly; this
includes all small prime factors, as well as $A$ itself if $A$ is prime (note that
rigorous primality proofs are not needed because correctness of the
parts of the algorithm that use primality can be verified afterwards).
Storing the $g$'s in previous steps may help; a partial
factorization of $A$ may already be known.
If any prime factors were found, store them, divide $A$ by them, and repeat.
\item Now $A$ no longer contains easily computable prime factors, which
implies that $A$ is not a prime.
Check if $A$ is a square. If so, store the square root.
Otherwise, apply a factorization algorithm to $A$ that finds at least one
non-trivial factor and go back to the previous step.
\item At this stage enough information has been computed to write $a = a_1 a_2^2$
where $a_1,a_2$ are integers, $a_1$ is square-free, and the factorization of $a_1$
(but in general not $a_2$) is known.
\item Apply the same process to $B$ and $C$.
\end{enumerate}
\item If $k>0$ then apply square-free factorization to write $a$ as $a_1 a_2^2$
for some $a_1, a_2 \in \tilde{R}_k$ with $a_1$ square-free.
Now factor $a_1$ in $\tilde{R}_k$.
Note that the {\em contents}\, of $a_1$ (the largest $r \in 
\Q$ for which $a_1/r \in R_k$) will usually be an integer. This integer does
not need to be square-free and should not be factored. Apply the same process to $b$ and $c$.
\item Return $a_1$, $b_1$, $c_1$ and stop.
\end{enumerate}
If $k=0$, each time a prime $p$ dividing $A$ is found, we can check
if $f_a$ can be factored mod $p$, since if it remains irreducible, then the conic
has no points over $\Q$ and so there is no reason to continue.
Apart from the output, one also needs to store $a_2,b_2,c_2$ and the $g$'s from
steps 2, 3, 4 since this information is needed to map solutions of the reduced equation
to solutions of the original equation. \\

\noindent {\bf Algorithm Conic} \\[2pt]
{\bf Input:} $k$, \ $t_1,\ldots,t_k$, \ \ $a,b,c \in K$. \\[2pt]
{\bf Output:} A solution $(X:Y:Z) \in P^2(K)$ of equation~(\ref{conic}) if such solution exists.
\begin{enumerate}
\item
If $a=0$ then return $(1:0:0)$ and stop. \\
If \hspace{1pt}$b=0$\hspace{0pt} then return $(0:1:0)$ and stop. \\
If \hspace{1pt}$c=0$\hspace{0pt} then return $(0:0:1)$ and stop.
\item If $k=0$ and $a,b,c$ have the same sign then return ``no solution exists'' and stop.
\item Replace $a,b,c$ by the output of Algorithm ReduceConic, and determine from the stored information
a map $\phi$ that maps solutions of the reduced equation to solutions of the original equation.
\item Compute a solubility certificate:
\begin{enumerate}
\item
For each $p$ in the factors of $a$ resp. $b$ resp. $c$,
factor the image of $f_a$ resp. $f_b$ resp. $f_c$ in ${\bf F}_p[T]$. If irreducible, then
return ``no solution exists'' and stop. Otherwise, store one of the roots (this root
will be used in step~\ref{step6} of Algorithm FindPoint).
\item
If ${\rm case}=0$ then, instead of computing a solubility certificate for equation~(\ref{lc}),
apply Algorithm Conic to find a solution $(x:y:z)$ of~(\ref{lc}).
If this recursive call returned ``no solution exists'' then return the same message and stop.
\end{enumerate}
\item Let $(X:Y:Z)$ be the output of Algorithm FindPoint.
\item Return $\phi(X:Y:Z)$ and stop.
\end{enumerate}

If $F$ is a field for which the necessary ingredients are available
(an implementation to solve a conic over $F$,
and the factorization algorithms mentioned in section~\ref{sect21})
then $F(t_1,\ldots,t_k)$ can be handled as well.
If a conic implementation over $F$ is not available,
Algorithm Conic in this paper could still be useful since
it would still work when ${\rm case}=1, \,k>0$.

\section{Examples}
An implementation in Maple is available at \\[1pt]
\verb+http://www.math.fsu.edu/~hoeij/files/ConicProgram+ \\[10pt]
{\bf Example 1.} \\[3pt]
Let $a := 1$, \ $b := 2t_2^2+ 2t_1 t_2-t_1^2$, \ and \ 
$c := -3t_2^4-4t_1 t_2^3  + 8 t_1^2 t_2^2 + 16 t_1^3 t_2 - 48t_1^4$. \\[3pt]
The implementation does not always find the same solution, see the remarks at the end
of section~\ref{sect21}. As output, one encounters for example
\[ (t_2^3+t_1 t_2^2 + 4t_1^2 t_2 - 4t_1^3 \ : \ t_2^2 - 4t_1^2 \ : \ t_2) \]
but also
\[ (t_2^3 + 2 t_1^2 t_2 + 8t_1^3 \ : \ t_2^2 + 2 t_1 t_2 -4t_1^2 \ : \ t_2+t_1). \]
\vspace{10pt}

\noindent {\bf Example 2.} \\[3pt]
$a := t_1^2+1$, \ $b := -(t_1^2+1)t_2^2 + (4t_1^2+4)t_2-1$, \ $c := 
(2t_1^3-10t_1^2+2t_1-9)t_2^4$ \\[3pt] $-(6t_1^3 + 6t_1-2)t_2^3+
(t_1^4-8t_1^3-2t_1^2-6t_1-2)t_2^2-(4t_1^4+2t_1^2-2t_1-2)t_2-1$. \\[3pt]
The implementation may find
\[ (t_2^3-t_2^2+t_2-1 \ : \ t_2^2+t_2+t_1 \ : \ 1) \]
but may also find a more complicated solution, depending on
which solubility certificate it chose.

\end{document}